\documentclass[a4paper,12pt]{article}

\usepackage{hyperref}
\usepackage{graphicx}
\usepackage{url}
\usepackage{amsthm}
\usepackage{newlfont}
\usepackage{setspace}
\usepackage{longtable}
\usepackage{ltxtable}
\usepackage{multirow}
\usepackage{amsmath}
\usepackage{amsfonts}
\usepackage{amssymb}

\newtheorem{theo}{Theorem}
\newtheorem{prop}{Proposition}
\newtheorem{lem}{Lemma}

\numberwithin{theo}{section}

\title{Generalization of Some Arithmetical Properties of Fermat-Euler Dynamical Systems}
\author{Ahmed N. Elsawy\footnote{Adresse: Mathematisches Institut,
Heinrich-Heine-Universit\"at,
D-40225 D\"usseldorf,
Germany, Email: elsawy@math.uni-duesseldorf.de}}

\begin{document}

\maketitle

\begin{abstract}
We study and generalize some arithmetical properties of the classes $(2^k+)$ and $(2^k-)$ introduced by V.~I.~Arnold:
a number $n$ belongs to the class $(N+)$ if $N | \varphi(n)$ and $2^{\frac{\varphi(n)}{N}} \equiv 1 \textrm{ mod } n$ where $\varphi(n)$ is 
the Euler function, and belongs to the class $(M-)$ if $M | \varphi(n)$ and $2^{\frac{\varphi(n)}{M}} \equiv -1 \textrm{ mod } n$. 
The classes $(2\pm),(4\pm)$ and $(8\pm)$ are studied by V.~I.~Arnold in \cite{Ferm-Eul,Ergodic-arithmatic} 
and here we will show general properties of the classes $(2^k\pm)$ and we will
see that the properties which proved by V.~I.~Arnold are special cases of ours.  
\end{abstract}

\section{Introduction}
One of the most important theorems in number theory and theory of finite groups is Fermat's little theorem 
($a^{p-1} \equiv 1 \textrm{ mod } p$, where $p$ is prime and $a \not \equiv 0\mod p$ is an integer) and 
its Euler's generalization ($a^{\varphi(n)} \equiv 1 \textrm{ mod } n$, where $a$ and $n$ are coprime integers).
They are also some of the most important applications of the dynamical system theory. The Fermat-Euler dynamical system is an example
of the connection between dynamical systems and number theory in which we try to use the properties of the dynamical systems to 
generalize or get new properties in number theory. 
\\[10pt]
\noindent {\bf Definition:}
The Abelian multiplicative group of residues coprime to $n$ is called the {\it Euler group} $\Gamma(n)$, and 
the {\it Euler function} $\varphi(n)$ is its order.

The Euler function is calculated as follows: 
 $$\varphi(1)=1 \textrm{  and } \varphi(n)= \prod (p_i-1)p^{a_i-1}_i,$$ where $n=\prod p^{a_i}_i$ 
is the prime factorization of $n$, and has the property
$\varphi(mn)\varphi(d)=\varphi(m)\varphi(n)d \textrm{  where } d=(m,n)$.
\\[10pt]
\noindent {\bf Definition:}
A Fermat-Euler dynamical system is a dynamical system with the ring of residues modulo $n$ as its set and the function $f(x)=ax$, 
where $(a,n)=1$, as its function.
\\[10pt]
\indent In what follows, we consider the case in which $a=2$ and $n$ is odd.
In \cite{Ferm-Eul}, V.~I.~Arnold has proved the following theorem which is a generalization of Euler's Theorem, by using 
the Fermat-Euler dynamical system.
\\[10pt]
\noindent {\bf Theorem A}\cite{Ferm-Eul}. All cycles of the Fermat-Euler permutation have the same length (period) $T$.
Therefore, the period and the number $N=N(n)$ of cycles satisfy the relation 
$$\varphi(n)=NT \quad \textrm{ and } \quad a^{\varphi(n)/N}\equiv 1 \mod n.$$ 

This theorem was the starting point of defining and studying new classes of numbers called $(N\pm)$.
\\[10pt]
\noindent {\bf Definition:} An odd number $n$  belongs to the class $(N+)$ if $N | \varphi(n)$ and 
$$2^{\varphi(n)/N} \equiv +1 \textrm{ mod } n,$$ and belongs to the class $(M-)$ if $M | \varphi(n)$ and 
$$2^{\varphi(n)/M} \equiv -1 \textrm{ mod } n.$$
The following properties directly follow from the definition 
\begin{enumerate}
 \item If $K \mid N$, then $(N+)\subseteq (K+)$.
 \item If $M$ is odd, then $(M-)$ is empty.
 \item If $M=K(2m+1)$, then $(M-)\subseteq (K-)$.
 \item $(2N-)\subseteq (N+)$.
\end{enumerate}
The proof of these properties can be found in \cite{Ferm-Eul}.  
\\[10pt]
\indent In \cite{Ferm-Eul}, V.~I.~Arnold calculated for each odd integer $1 < n < 512$ its classes $(N+)$ and $(M-)$ with
maximal values of $N$ and $M$ as well as the values of its minimal period $T(n)$ and the Euler function $\varphi(n)$.
From these calculations he observed and proved some properties of the classes $(2\pm),(4\pm)$ and $(8\pm)$. 
In this work we give and prove general properties of the classes $(2^k+)$ and $(2^k-)$, where $k$ is a positive integer, and show 
that the properties which are proved in \cite{Ferm-Eul,Ergodic-arithmatic} by V.~I.~Arnold are special cases of the generalized properties.
\\[10pt]
\indent In \cite{Ergodic-arithmatic}, V.~I.~Arnold proved by induction that if $n$ is divisible by more than $k$ distinct primes, then $n\in (2^k+)$.
We will give here another proof and we will answer the following question. 
\\[10pt]
\noindent {\bf Question:}
When does $n$ belong to $(2^k+)$ or $(2^k-)$ if $n$ has exactly $k$ or $k-1$ distinct prime divisors, and when does it not?
\\[10pt]
In section 2, we prove the following theorems which give the complete answer to this question.
\\[10pt]
\noindent {\bf Theorem 1.}
Let $n$ be an odd integer and have exactly $k\geq 2$ distinct prime divisors.
\begin{itemize}
\item If $n=\prod^k_{i=1} p^{a_i}_i$ with $p_i\equiv 3\mod 8$ for all $i$, then $n\in (2^k-)$.

\item If $n=\prod p^{a_i}_i\prod q^{b_j}_j$ with $p_i\equiv 3 \mod 8$ for all $i$, and $q_j\equiv -1\mod 8$ for all $j$,
then $ n \notin (2^k+)\cup(2^k-)$.

\item If $n=q^b\prod_{i=1}^{k-1} p^{a_i}_i$ with $p_i\equiv 3 \textrm{ or }-1\mod 8$, and 
$q\equiv -3\mod 8$, then again $n \notin (2^k+)\cup(2^k-)$.

\item If $n$ is not in the previous three cases, then $n \in (2^k+)$.
\end{itemize}

\noindent {\bf Theorem 2.}
Let $n$ be an odd integer and have exactly $k-1\geq 2$ distinct prime divisors.
\begin{itemize}
\item If $k=3$ and $n=p_1^{a_1}p_2^{a_2}$ with $p_1\equiv p_2\equiv -3\mod 8$, then $n\in (2^k-)$.

\item If $n=\prod^{r}_{i=1}p_i^{a_i}\prod^{k-1}_{j=r+1}q_j^{b_j}$ with $p_i\equiv -3 \mod 8$ for all $i$,
$q_j\equiv -1 \textrm{ or } 3 \mod 8$ for all $j$, and $0\leq r \leq 2$, then
$n\notin (2^k+)\cup(2^k-)$.

\item $n=p_1^{a_1}\prod^{k-1}_{j=2}q_j^{b_j}$ with 
$p_1\in (4-)$ and $q_j\equiv -1 \textrm{ or } 3 \mod 8$ for all $j$, then $n\notin (2^k+)\cup(2^k-)$.

\item If $n$ is not in the previous three cases, then $n \in (2^k+)$.
\end{itemize}

In section 3 we will see the cases $k=2$ and $k=3$ as an example and how the theorems in section 2 describe the classes
$(4\pm)$ and $(8\pm)$. Also we will compare between our results and the results which are proven by V.~I.~Arnold in 
\cite{Ferm-Eul,Ergodic-arithmatic}.

\section{The classes ($2^k\pm$)} 
\label{sec:2}

We will use the following propositions in our proofs, 
which can be found in any elementary book of Number Theory, for example in \cite{Friendly}.
\begin{prop}
\label{prop:1}
If $x\equiv a \mod A \textrm{ and } x\equiv a \mod B, \textrm{ then } x\equiv a \mod C$, where $C$ is the least common multiple of 
$A$ and $B$. 
\end{prop}

\begin{prop}
\label{prop:2}
Let $p$ be an odd prime, then
\begin{displaymath}
 2^{\frac{\varphi(p^a)}{2}} \equiv \left\{ \begin{array}{ll}
                          +1 \mod p^a & \textrm{if } p \equiv \pm1 \mod 8 \\
                          -1 \mod p^a & \textrm{if } p \equiv \pm3 \mod 8
                         \end{array}\right.
\end{displaymath} 
\end{prop}

This proposition can be easily obtained from Euler's criterion 
($a^{\frac{p-1}{2}} \equiv (\frac{a}{p}) \mod p, \textrm{ where }(\frac{a}{p})$ is the Legendre symbol). 
Indeed, since we know that
\begin{displaymath}
 \Big (\frac{2}{p} \Big ) \equiv \left\{ \begin{array}{ll}
                          +1 \mod p & \textrm{if } p \equiv \pm1 \mod 8 \\
                          -1 \mod p & \textrm{if } p \equiv \pm3 \mod 8,
                         \end{array}\right.
\end{displaymath}
we can write $2^{\frac{p-1}{2}}=\pm1+Ap$ for some integer $A$. 
Hence, using binomial formula, we can easily get 
$$2^{\frac{(p-1)p^{a-1}}{2}}=(\pm1+Ap)^{p^{a-1}}=\pm1+A'p^a,\textrm{ because } p^{a-1} \textrm{ is odd }.$$

The following theorem is proved in \cite{Ergodic-arithmatic} by induction, we give here a direct proof
using Proposition 2.
\begin{theo}
\label{th:1}
 Every odd number which is divisible by more than $k$ different primes belongs to the class $(2^k+)$.
\end{theo}
 \proof
Let $n$ be an odd number which is divisible by exactly $r$ different primes with $0<k<r$, then we can write $n$ as
$n=AB$ where $A=\prod^k_{i=1} p^{a_i}_i$ and $B=\prod^r_{i=k+1}p^{a_i}_i$. 

From Euler's theorem, it follows that $2^{\varphi(B)}\equiv 1 \mod B$. 
Also, note that $\varphi(A)/2^k$ is an integer, because $A$ has $k$ different (odd) prime divisors.
So 
$$2^{\varphi(n)/2^k}=2^{\varphi(A)\varphi(B)/2^k}=(2^{\varphi(B)})^{\varphi(A)/2^k}\equiv 1 \mod B.$$

Since $\varphi(B)$ is even, $ \varphi(A)=\prod_{i=1}^k\varphi(p^{a_i}_i)$, and $2^{\varphi(p^{a_i}_i)/2} \equiv \pm1 \mod p^{a_i}_i$,
we have
\begin{displaymath}
2^{\varphi(n)/2^k}=2^{\varphi(A)\varphi(B)/2^k}=
(2^{\varphi(p^{a_i}_i)/2})^{\varphi(B)\frac{\prod_{j=1, j\neq i}^k\varphi(p^{a_j}_j)}{2^{k-1}}}\equiv 1\mod p^{a_i}_i
\end{displaymath}
$$\textrm{ for every } 1\leq i \leq k.$$
Thus, using Proposition 1, we get $2^{\varphi(n)/2^k}\equiv 1\mod \prod_{i=1}^k p^{a_i}_i=A$.

Now we have $2^{\varphi(n)/2^k}\equiv 1\mod A$ and $2^{\varphi(n)/2^k}\equiv 1\mod B$ where $A$ and $B$ are relatively prime; 
therefore,
$$\qquad\qquad\quad 2^{\varphi(n)/2^k}\equiv +1\mod AB=n \textrm{; that is, } n \in (2^k+).\qquad\qquad\qquad\Box$$ 

Theorem \ref{th:1} and Proposition 2 give us, in particular, a complete description of the classes $(2\pm)$, i.e. the case $k=1$, as follows.
\begin{itemize}
\item Every odd number which is divisible by more than one prime belongs to the class $(2+)$.

\item If $n=p^a$ with $p\equiv \pm1$ mod 8, then $n\in (2+)$.

\item If $n=p^a$ with $p\equiv \pm3$ mod 8, then $n\in (2-)$.
\end{itemize}

Now the question is, what will happen in the case when $n$ has exactly $k$ distinct prime divisors? 
In Theorems \ref{th:2}-\ref{th:6} we study this case and show  precisely when $n$ belongs to the class $(2^k+)$, to the class $(2^k-)$, and 
when $n$ belongs to neither of them.

\begin{theo}
\label{th:2}
If $n=\prod^k_{i=1} p^{a_i}_i$ has at least one prime divisor congruent to 1 modulo 8, then $n\in (2^k+)$.
\end{theo}
\proof
WLOG, we assume that $p_1\equiv 1 \mod 8$ and $n=p_1^{a_1}A$ where $A=\prod^k_{i=2} p^{a_i}_i$. Then

$$2^{\varphi(n)/2^k}=(2^{\frac{\varphi(p_1^{a_1})}{2}})^{\frac{\varphi(A)}{2^{k-1}}}\equiv (+1)^{\frac{\varphi(A)}{2^{k-1}}}\equiv 1 \mod p_1^{a_1}.$$
For all $i=2,\dots, k$, we have $2^{\varphi(p^{a_i}_i)/2}\equiv +1\textrm{ or }-1\mod p^{a_i}_i$, and 
since $p_1\equiv 1\mod 8$, it follows that $\frac{p_1-1}{2}$ is even, so $\frac{\varphi(p_1^{a_1})}{2}$ is even. Thus,
$$2^{\varphi(n)/2^k}=(2^{\frac{\varphi(A)}{2^{k-1}}})^{\frac{\varphi(p_1^{a_1})}{2}}=
(2^{\varphi(p^{a_i}_i)/2})^{\frac{\varphi(p_1^{a_1})}{2}\frac{\prod_{j=2, j\neq i}^k\varphi(p^{a_j}_j)}{2^{k-2}}}\equiv 1\mod p^{a_i}_i$$
$$\textrm{ for every }  2\leq i \leq k.$$
It follows that $2^{\varphi(n)/2^k}\equiv +1 \mod A$, by Proposition 1.
\\[10pt]
So we have $2^{\varphi(n)/2^k}\equiv 1 \mod p_1^{a_1}$ and $2^{\varphi(n)/2^k}\equiv 1 \mod A$ 
where $p_1^{a_1}$ and $A$ are relatively prime.
Thus, 

$2^{\varphi(n)/2^k}\equiv +1 \mod p_1^{a_1}A=n \textrm{; that is,  }n \in (2^k+)$. \hspace{\stretch{1}} $\Box$
\\[10pt]
\indent Now assume that $n$ has no prime divisors congruent to 1 modulo 8. 
In Theorems \ref{th:3}, and \ref{th:4}, we deal with the case when all the prime divisors of $n$ have the same 
residue (except 1) modulo 8. 
\begin{theo}
\label{th:3}
If $n=\prod^k_{i=1} p^{a_i}_i$ with either $p_i\equiv -1\mod 8$ for all $i$ or $p_i\equiv -3\mod 8$ for all $i$, then $n\in (2^k+)$.
\end{theo}
\proof
Case 1. Let $p_i\equiv -1\mod 8$ for all $i$, then by Proposition 2 we have 
$2^{\varphi(p^{a_i}_i)/2}\equiv 1 \mod p^{a_i}_i$ for all $i$, so
$$2^{\varphi(n)/2^k}=(2^{\varphi(p^{a_i}_i)/2})^{\frac{\prod_{j=1, j\neq i}^k\varphi(p^{a_j}_j)}{2^{k-1}}}\equiv 
(+1)^{\frac{\prod_{j=1, j\neq i}^k\varphi(p^{a_j}_j)}{2^{k-1}}}\equiv 1 \mod p^{a_i}_i.$$
Thus, $2^{\varphi(n)/2^k}\equiv +1 \mod n$; that is, $n \in (2^k+)$.
\\[10pt]
Case 2. Let $p_i\equiv -3\mod 8$ for all $i$, then $2^{\varphi(p^{a_i}_i)/2}\equiv -1 \mod p^{a_i}_i$, and
$\varphi(p^{a_i}_i)/2$ is even, for all $i$. Indeed, let $p_i=5+r_i8$ for some integer $r_i$, then $\frac{p_i-1}{2}=2+r_i4$ is even;
therefore, $\varphi(p_i^{a_i})/2=p_i^{a_i-1}(p_i-1)/2$ is even for all $i$.

Hence,$$2^{\varphi(n)/2^k}=(2^{\varphi(p^{a_i}_i)/2})^{\frac{\prod_{j=1, j\neq i}^k\varphi(p^{a_j}_j)}{2^{k-1}}}\equiv 
(-1)^{\frac{\prod_{j=1, j\neq i}^k\varphi(p^{a_j}_j)}{2^{k-1}}}\equiv +1 \mod p^{a_i}_i.$$
So we have again $2^{\varphi(n)/2^k}\equiv +1 \mod n$, or $n \in (2^k+)$.\hspace{\stretch{1}} $\Box$

\begin{theo}
\label{th:4}
If $n=\prod^k_{i=1} p^{a_i}_i$ with $p_i\equiv 3\mod 8$ for all $i$, then $n\in (2^k-)$.
\end{theo}
\proof
Since $p_i\equiv 3\mod 8$ for all $i$, it follows from Proposition 2 that $2^{\varphi(p^{a_i}_i)/2}\equiv -1 \mod p^{a_i}_i$, and
$\varphi(p^{a_i}_i)/2$ is odd, for all $i$. Indeed, let $p_i=3+r_i8$ for some integer $r_i$, then $\frac{p_i-1}{2}=1+r_i4$ is odd; therefore,
$\varphi(p_i^{a_i})/2=p_i^{a_i-1}(p_i-1)/2$ is odd for all $i$.

Thus,
$$2^{\varphi(n)/2^k}=(2^{\varphi(p^{a_i}_i)/2})^{\frac{\prod_{j=1, j\neq i}^k\varphi(p^{a_j}_j)}{2^{k-1}}}\equiv 
(-1)^{\frac{\prod_{j=1, j\neq i}^k\varphi(p^{a_j}_j)}{2^{k-1}}}\equiv -1 \mod p^{a_i}_i.$$
So $2^{\varphi(n)/2^k}\equiv -1 \mod \prod^k_{i=1} p^{a_i}_i=n$; that is, $n \in (2^k-)$.\hspace{\stretch{1}} $\Box$
\\[10pt]
In the next two theorems, $n$ has no prime divisor $p\equiv 1 \mod 8$ and has at least two prime divisors $p$ and $q$ with 
different residues modulo 8.
\begin{theo}
\label{th:5}
Let $n=\prod_{i=1}^r p^{a_i}_i\prod_{j=r+1}^k q^{b_j}_j$ with $ 0<r<k,\ p_i\equiv 3$ or $-1$ mod 8 for all $i$, and 
$q_j\equiv -3$ mod 8 for all $j$, then 
$$ \begin {array}{ll}
n \notin (2^k+)\cup(2^k-)& \textrm{  if  } r=k-1, \textrm{ and }\\    
n \in (2^k+)             &  \textrm{  if  } r<k-1.
\end{array}$$
\end{theo}
\proof
Note that $\varphi(p_i^{a_i})/2$ is odd. Indeed, $p_i\equiv$ 3 or $-1$ mod 8 means that $p_i\equiv 3 \mod 4$, so we can write
$p_i=3+4s_i$ for some integer $s_i$; therefore, $\frac{p_i-1}{2}=1+2s_i$, and $\varphi(p_i^{a_i})/2=p_i^{a_i-1}(\frac{p_i-1}{2})$ is odd.
\\[10pt]
On the other hand, $\varphi(q_j^{b_j})/2$ is even. Indeed, $q_j\equiv -3$ mod 8 means that $q_j\equiv 1 \mod 4$, and we can write
$q_j=1+4s_j$ for some integer $s_j$; therefore, $\frac{q_j-1}{2}=2s_j$, so $\varphi(q_j^{b_j})/2=q_j^{b_j-1}(\frac{q_j-1}{2})$ is even.
\\[10pt]
So we have
$$2^{\varphi(n)/2^k}=(2^{\varphi(p_i^{a_i})/2})^{\frac{\prod_{l=1,l\neq i}^r \varphi(p^{a_l}_l)\prod_{j=r+1}^k
\varphi(q^{b_j}_j)}{2^{k-1}}}\equiv (\pm1)^{\frac{\prod_{l=1,l\neq i}^r \varphi(p^{a_l}_l)\prod_{j=r+1}^k
\varphi(q^{b_j}_j)}{2^{k-1}}}$$
$$\equiv +1 \mod p_i^{a_i}, \textrm{  because } \varphi(q_j^{b_j})/2 \textrm{ is even}.$$
Therefore, $2^{\varphi(n)/2^k}\equiv 1 \mod \prod_{i=1}^r p_i^{a_i}$.
\\[10pt]
Case 1. Let $r=k-1$; that is, $n=q^b\prod_{i=1}^{k-1} p^{a_i}_i$ with $p_i\equiv 3$ or $-1$ mod 8, and 
$q\equiv -3\mod 8$, then 
$$2^{\varphi(n)/2^k}=(2^{\varphi(q^b)/2})^{\frac{\prod_{i=1}^{k-1}\varphi(p^{a_i}_i)}{2^{k-1}}}\equiv 
(-1)^{\frac{\prod_{i=1}^{k-1}\varphi(p^{a_i}_i)}{2^{k-1}}}\equiv -1 \mod q^b,$$
because $\frac{\prod_{i=1}^{k-1}\varphi(p^{a_i}_i)}{2^{k-1}}$ is odd. 
Thus, 

$2^{\varphi(n)/2^k}\equiv +1 \mod \prod_{i=1}^r p_i^{a_i}$ and $2^{\varphi(n)/2^k}\equiv -1 \mod q^b$.

\noindent Now assume that $n\in (2^k+)\cup(2^k-)$, then $2^{\varphi(n)/2^k}\equiv \pm1 \mod n$;
therefore, 
$2^{\varphi(n)/2^k}$ has the same residue modulo $\prod_{i=1}^r p_i^{a_i}$ 
and modulo $q^b$. But this contradicts our calculations, so $n\notin (2^k+)\cup(2^k-)$.
\\[10pt]
Case 2. Let $0<r<k-1$ and $n=\prod_{i=1}^r p^{a_i}_i\prod_{j=r+1}^k q^{b_j}_j$, then 
$$2^{\varphi(n)/2^k}=(2^{\varphi(q_j^{b_j})/2})^{\frac{\prod_{i=1}^r\varphi(p^{a_i}_i)\prod_{l=r+1,l\neq j}^k\varphi(q^{b_l}_l)}{2^{k-1}}}
\equiv (-1)^{\frac{\prod_{i=1}^r\varphi(p^{a_i}_i)\prod_{l=r+1,l\neq j}^k\varphi(q^{b_l}_l)}{2^{k-1}}}$$
$$\equiv +1 \mod q_j^{b_j},$$
because $\frac{\varphi(q^{b_l}_l)}{2}$ is even.
Thus, 

$2^{\varphi(n)/2^k}\equiv 1 \mod \prod_{i=1}^r p_i^{a_i}$ and $2^{\varphi(n)/2^k}\equiv 1 \mod \prod_{j=r+1}^k q^{b_j}_j$,
which means that 
$2^{\varphi(n)/2^k}\equiv +1 \mod \prod_{i=1}^r p_i^{a_i}\prod_{j=r+1}^k q^{b_j}_j=n$, or $n\in (2^k+)$.\hspace{\stretch{1}} $\Box$
\\[10pt]
Note that if $r$ in Theorem \ref{th:5} was equal to zero, then we would have Theorem \ref{th:3}. So the last possibility is that 
$r$ is equal to $k$, which we investigate in the following theorem.

\begin{theo}
\label{th:6}
Let $n=\prod p^{a_i}_i \prod q^{b_j}_j$ with $p_i\equiv 3 \mod 8$ for all $i$, and $q_j\equiv -1\mod 8$ for all $j$, then 
$$ n \notin (2^k+)\cup(2^k-).$$
\end{theo}
\proof
Both $\varphi(p_i^{a_i})/2$ and $\varphi(q_j^{b_j})/2$ are odd, because $p_i\equiv q_j\equiv 3 \mod 4$, and so
$\frac{p_i-1}{2}\equiv \frac{q_j-1}{2}\equiv 1 \mod 2$. So we have
$$2^{\varphi(n)/2^k}=(2^{\varphi(p_i^{a_i})/2})^{\frac{\prod_{l\neq i}\varphi(p^{a_l}_l)\prod\varphi(q^{b_j}_j)}{2^{k-1}}}\equiv 
(-1)^{\frac{\prod_{l\neq i}\varphi(p^{a_l}_l)\prod\varphi(q^{b_j}_j)}{2^{k-1}}}\equiv -1 \mod p_i^{a_i},$$
because $\frac{\prod_{l\neq i}\varphi(p^{a_l}_l)\prod\varphi(q^{b_j}_j)}{2^{k-1}}$ is odd.
On the other hand,
$$2^{\varphi(n)/2^k}=(2^{\varphi(q_j^{b_j})/2})^{\frac{\prod_{l\neq j}\varphi(p^{a_i}_i)\prod\varphi(q^{b_l}_l)}{2^{k-1}}}\equiv 
(+1)^{\frac{\prod_{l\neq j}\varphi(p^{a_i}_i)\prod\varphi(q^{b_l}_l)}{2^{k-1}}}\equiv 1 \mod q_j^{b_j}.$$
Thus, $2^{\varphi(n)/2^k}\equiv -1 \mod \prod p_i^{a_i}$ and $2^{\varphi(n)/2^k}\equiv +1 \mod \prod q_j^{b_j}$, 
which contradicts both $2^{\varphi(n)/2^k}\equiv +1 \mod n$ and $2^{\varphi(n)/2^k}\equiv -1 \mod n$; that is, 
$n \notin (2^k+)\cup(2^k-).$ \hspace{\stretch{1}} $\Box$
\\[10pt]
Note that Theorems \ref{th:2}, \ref{th:3}, \ref{th:4}, \ref{th:5}, and \ref{th:6} 
cover all possible cases of $n$ with $k\geq 2$ different prime divisors, and so Theorem 1 in the introduction is proved.
\\[10pt]
So the case when $n$ has exactly $k$ prime divisors is completely described. What can we say about the case when $n$ has exactly
$k-1$ prime divisors? Does $n$ belong or not belong to the classes $(2^k+)$ or
$(2^k-)$ in this case?
Theorems \ref{th:7}, \ref{th:8}, \ref{th:9}, and \ref{th:10} give a complete answer to this question.

From here on, we consider the case
$$n=\prod^{k-1}_{i=1}p_i^{a_i} \textrm{ where } k\geq 3.$$ 
In Theorems \ref{th:7}, and \ref{th:8}, we study the case in which $n$ has at least one prime divisor congruent to 1 modulo 8.
In Theorem \ref{th:7}, we assume that $n$ has another prime divisor congruent to either 1 or $-3$ modulo 8, and 
in Theorem \ref{th:8}, we assume that all prime divisors but one, which is congruent to 1 as we already mentioned, 
are congruent to 1 or $-3$ modulo 8.
\begin{theo}
\label{th:7}
If $n$ has at least two prime divisors, say  $p_1$ and $p_2$, with $p_1\equiv 1\mod 8$ and $p_2\equiv$ 1 or $-3$ mod 8, then 
$n \in (2^k+)$. 
\end{theo}
\proof
Let $n=p_1^{a_1}p_2^{a_2}A$, where $A=\prod^{k-1}_{i=3}p_i^{a_i}$  and  $k\geq 3$.
Since $p_1\equiv 1\mod 8$ and moreover, $\frac{\varphi(p_2^{a_2})}{4}$ is an integer (since $p_2\equiv$ 1 or $-3$ mod 8 means that 
$p_2-1$ is divisible by 4),
we have
$$2^{\frac{\varphi(n)}{2^k}}=(2^{\frac{\varphi(p_1^{a_1})}{2}})^{\frac{\varphi(p_2^{a_2})}{4}\frac{\varphi(A)}{2^{k-3}}}
\equiv 1 \mod p_1^{a_1}.$$
Since $\frac{\varphi(p_1^{a_1})}{4}$ is even ($8|p_1-1$), we have
$$2^{\frac{\varphi(n)}{2^k}}=
(2^{\frac{\varphi(p_i^{a_i})}{2}})^{\frac{\varphi(p_1^{a_1})}{4}\frac{\prod^{k-1}_{j=2,j\neq i}\varphi(p_j^{a_j})}{2^{k-3}}}
\equiv 1 \mod p_i^{a_i} \textrm{ for all } i=2,\dots, k-1.$$
Thus, 
$2^{\frac{\varphi(n)}{2^k}}\equiv 1 \mod p_i^{a_i} \textrm{ for all } i=1,\dots, k-1$; therefore,  
$2^{\frac{\varphi(n)}{2^k}}\equiv +1 \mod n$, in other words, $n \in (2^k+)$. \hspace{\stretch{1}} $\Box$
\\[10pt]
The following lemma give us an important property of the primes in the classes $(4+)$ and $(4-)$, 
which we need in the proof of Theorem \ref{th:8}.
\begin{lem}
\label{th:3.1}
A prime power $p^a \in (4+)\cup(4-)$ if and only if $p \equiv 1\mod 8$.
\end{lem}
\proof
First, if $p^a \in (4+)\cup(4-)$, then $p \equiv 1\mod 8$, is already proved in \cite{Ferm-Eul}.
Conversely, if $p \equiv 1\mod 8$, then $2^{\frac{p-1}{2}}\equiv 1 \mod p$. So we can write
$$0\equiv 2^{\frac{p-1}{2}}-1\equiv (2^{\frac{p-1}{4}}-1)(2^{\frac{p-1}{4}}+1) \mod p.$$ 
Which means that $2^{\frac{p-1}{4}}\equiv +1 \textrm{ or } -1 \mod p$, that is $2^{\frac{p-1}{4}}=\pm1+Ap$ for some integer $A$.
Thus
$$2^{\frac{\varphi(p^a)}{4}}=(2^{\frac{p-1}{4}})^{p^{a-1}}=(\pm1+Ap)^{p^{a-1}}=\pm1+Bp^a \textrm{ for some integer } B.$$
Then $p^a \in (4+)\cup(4-).$ \hspace{\stretch{1}} $\Box$

\begin{theo}
\label{th:8}
If $n=p_1^{a_1}\prod^{k-1}_{i=2}p_i^{a_i}$ with $p_1\equiv 1\mod 8$ and $p_i\equiv -1$ or $3$ 

\noindent mod 8 for all $2\leq i\leq k-1$, then 
$$\begin{array}{ll}
n \in (2^k+) & \textrm{ if } p_1\in (4+), \textrm{ and } \\
n \notin (2^k+)\cup(2^k-) & \textrm{ if } p_1\in (4-).
  \end{array}
$$ 
\end{theo}
\proof
Let $n=p_1^{a_1}A$ where $A=\prod^{k-1}_{i=2}p_i^{a_i}$ with $p_i\equiv -1$ or 3 mod 8.

\noindent If $p_i\equiv -1$ mod 8, then $2^{\frac{\varphi(p_i^{a_i})}{2}}\equiv 1\mod p_i^{a_i}$. Thus,
$$2^{\frac{\varphi(n)}{2^k}}=
(2^{\frac{\varphi(p_i^{a_i})}{2}})^{\frac{\varphi(p_1^{a_1})}{4}\frac{\prod^{k-1}_{j=2,j\neq i}\varphi(p_j^{a_j})}{2^{k-3}}}
\equiv 1 \mod p_i^{a_i}.$$
If $p_i\equiv$ 3 mod 8, then $2^{\frac{\varphi(p_i^{a_i})}{2}}\equiv -1\mod p_i^{a_i}$. Since $\frac{\varphi(p_1^{a_1})}{4}$ is
even, we have 
$$2^{\frac{\varphi(n)}{2^k}}=
(2^{\frac{\varphi(p_i^{a_i})}{2}})^{\frac{\varphi(p_1^{a_1})}{4}\frac{\prod^{k-1}_{j=2,j\neq i}\varphi(p_j^{a_j})}{2^{k-3}}}
\equiv 1 \mod p_i^{a_i}.$$
So $2^{\frac{\varphi(n)}{2^k}}\equiv 1 \mod p_i^{a_i} \textrm{ for all } i=2,\dots, k-1$. Therefore, 
$$2^{\frac{\varphi(n)}{2^k}}\equiv 1 \mod A.$$
Since $\frac{\varphi(p_i^{a_i})}{2}$ is odd for all $i=2,\dots, k-1$, only $\varphi(p_1^{a_1})$ can be divisible by 4. Since 
$p^a \in (4+)\cup(4-)$ if and only if $p \equiv 1\mod 8$, as we proved in Lemma \ref{th:3.1}, we have $p_1\in (4+)$
or $\in (4-)$. 

\noindent If $p_1\in (4+)$, then $2^{\frac{\varphi(p_1^{a_1})}{4}}\equiv +1 \mod p_1^{a_1}$. Therefore,
$$2^{\frac{\varphi(n)}{2^k}}=
(2^{\frac{\varphi(p_1^{a_1})}{4}})^{\frac{\varphi(A)}{2^{k-2}}}
\equiv (+1)^{\frac{\varphi(A)}{2^{k-2}}} \equiv 1 \mod p_1^{a_1}.$$
So now we have $2^{\frac{\varphi(n)}{2^k}}\equiv 1 \mod p_1^{a_1}$ and 
$2^{\frac{\varphi(n)}{2^k}}\equiv 1 \mod A$, which means that $2^{\frac{\varphi(n)}{2^k}}\equiv +1 \mod p_1^{a_1}A=n$,
or  $n \in (2^k+)$.
\\[10pt]
If $p_1\in (4-)$, then $2^{\frac{\varphi(p_1^{a_1})}{4}}\equiv -1 \mod p_1^{a_1}$. Therefore,
$$2^{\frac{\varphi(n)}{2^k}}=(2^{\frac{\varphi(p_1^{a_1})}{4}})^{\frac{\varphi(A)}{2^{k-2}}}
\equiv (-1)^{\frac{\varphi(A)}{2^{k-2}}} \equiv -1 \mod p_1^{a_1},$$ because $\frac{\varphi(A)}{2^{k-2}}$ is odd. 
So we have $2^{\frac{\varphi(n)}{2^k}}\equiv -1 \mod p_1^{a_1}$ and 
$2^{\frac{\varphi(n)}{2^k}}\equiv +1 \mod A$. So $n \notin (2^k+)\cup(2^k-)$. Indeed, from $n \in (2^k+)\cup(2^k-)$, 
it follows that $2^{\frac{\varphi(n)}{2^k}}\equiv \pm1 \mod p_1^{a_1}$ and 
$2^{\frac{\varphi(n)}{2^k}}\equiv \pm1 \mod A$, which contradicts what we have. 
\hspace{\stretch{1}} $\Box$
\\[10pt]
{\bf Remark:} Since $p\equiv 1$ mod 8 is equivalent to $p \in (4+)\cup(4-)$, we can reformulate Theorems \ref{th:7} and
 \ref{th:8} as follows.
\begin{itemize}
\item If $n$ has at least one prime divisor $p \in (4+)$, then $n \in (2^k+)$.

\item If $n$ has at least two prime divisors $p_1$ and $p_2$ with $p_1\in (4-)$ and $p_2\equiv$ 1 or $-3$ mod 8, then 
$n \in (2^k+)$.

\item If $n$ has one prime divisor $p_1$ with $p_1 \in (4-)$ and all the other prime divisors $p_i\equiv -1$ or 3 mod 8, then
$n \notin (2^k+)\cup(2^k-)$.
\end{itemize}
finally, in Theorems \ref{th:9} and \ref{th:10} we consider the case in which $n$ has no prime divisor congruent to 1 modulo 8.
\begin{theo}
\label{th:9}
If $n=\prod^{k-1}_{i=1}p_i^{a_i}$ with $p_i\equiv -3\mod 8$ for all $i$, then 
$$\begin{array}{ll}
n \in (2^k+) & \textrm{ if } k > 3, \textrm{ and }\\
n \in (2^k-) & \textrm{ if } k = 3.
 \end{array}
$$ 
\end{theo}
\proof
Since $p_i\equiv -3\mod 8$, it follows that $2^{\frac{\varphi(p_i^{a_i})}{2}}\equiv -1 \mod p_i^{a_i}$, 
$\frac{\varphi(p_i^{a_i})}{2}$ is even, but $\frac{\varphi(p_i^{a_i})}{4}$ is odd. Hence,

$$2^{\frac{\varphi(n)}{2^k}}=(2^{\frac{\varphi(p_i^{a_i})}{2}})^{\frac{\varphi(p_l^{a_l})}{4}\frac{\prod^{k-1}_{j=1,j\neq i, l}\varphi(p_j^{a_j})}{2^{k-3}}}
\equiv (-1)^{\frac{\varphi(p_l^{a_l})}{4}\frac{\prod^{k-1}_{j=1,j\neq i, l}\varphi(p_j^{a_j})}{2^{k-3}}} \mod p_i^{a_i}.$$
If $k > 3$, then $\frac{\prod^{k-1}_{j=1,j\neq i, l}\varphi(p_j^{a_j})}{2^{k-3}}$ is even.
Thus, $2^{\frac{\varphi(n)}{2^k}} \equiv +1 \mod p_i^{a_i}$ for all $i$, and then 
$2^{\frac{\varphi(n)}{2^k}} \equiv +1 \mod n$, or $n\in (2^k+)$.
\\[10pt]
If $k=3$, then $n=p_1^{a_1}p_2^{a_2}$. Therefore, 

$2^{\frac{\varphi(n)}{2^3}}=(2^{\frac{\varphi(p_1^{a_1})}{2}})^{\frac{\varphi(p_2^{a_2})}{4}}\equiv -1 \mod p_1^{a_1}$, and
$2^{\frac{\varphi(n)}{2^3}}=(2^{\frac{\varphi(p_2^{a_2})}{2}})^{\frac{\varphi(p_1^{a_1})}{4}}\equiv -1 \mod p_2^{a_2}$.

Thus, $2^{\frac{\varphi(n)}{2^3}}\equiv -1 \mod p_1^{a_1}p_2^{a_2}=n$, or $n\in (2^k-)$. \hspace{\stretch{1}} $\Box$
\begin{theo}
\label{th:10}
If $n=\prod^{r}_{i=1}p_i^{a_i}\prod^{k-1}_{j=r+1}q_j^{b_j}$ with $0\leq r < k-1,\ p_i\equiv -3 \mod 8$ for all $i$, and
$q_j\equiv -1 \textrm{ or } 3 \mod 8$ for all $j$, then
$$\begin{array}{ll}
   n \notin (2^k+)\cup(2^k-)   & \textrm{ if  } 0 \leq r \leq 2, \textrm{ and }\\
   n \in (2^k+)                & \textrm{ if  } 2 < r < k-1.
  \end{array}
$$ 
\end{theo}
\proof
Now we have that $\varphi(p_i^{a_i})/2$ is even and $\varphi(p_i^{a_i})/4$ is odd. Indeed, $p_i\equiv -3 \mod 8$ means that 
$p_i=5+8s_i$ for some integer $s_i$, it follows that $\frac{p_i-1}{2}=2+4s_i$ is even and $\frac{p_i-1}{4}=1+2s_i$ is odd. Therefore, 
$\frac{\varphi(p_i^{a_i})}{2}=\frac{p_i-1}{2}p_i^{a_i-1}$ is even and $\frac{\varphi(p_i^{a_i})}{4}=\frac{p_i-1}{4}p_i^{a_i-1}$ is odd.

Moreover, $\varphi(q_j^{b_j})/2$ is odd. Indeed, $q_j\equiv -1 \textrm{ or } 3 \mod 8$ means that $q_j\equiv 3 \mod 4$ for all $j$, so
$\frac{q_j-1}{2}$ is odd, and therefore $\frac{q_j-1}{2}q_j^{b_j-1}=\varphi(q_j^{b_j})/2$ is odd.

If $r=0$, then $n=\prod^{k-1}_{j=1}q_j^{b_j}$ and $\varphi(n)=\prod^{k-1}_{j=1}\varphi(q_j^{b_j})$. Since 
$\varphi(q_j^{b_j})/2$ is odd for all $j$, it follows that $\frac{\prod^{k-1}_{j=1}\varphi(q_j^{b_j})}{2^{k-1}}$ is odd.
Thus, $\varphi(n)$ is not divisible by $2^k$, which implies that $n \notin (2^k+)\cup(2^k-)$.

If $r=1$, then $n=p_1^{a_1}\prod^{k-1}_{j=2}q_j^{b_j}$ and $\frac{\varphi(n)}{2^k}=
\frac{\varphi(p_1^{a_1})}{4}\frac{\prod^{k-1}_{j=2}\varphi(q_j^{b_j})}{2^{k-2}}$. So
$$2^{\frac{\varphi(n)}{2^k}}=(2^{\frac{\varphi(q_j^{b_j})}{2}})^
{\frac{\varphi(p_1^{a_1})}{4}\frac{\prod^{k-1}_{l=2,l\neq j}\varphi(q_l^{b_l})}{2^{k-3}}}\equiv (\pm 1)^
{\frac{\varphi(p_1^{a_1})}{4}\frac{\prod^{k-1}_{l=2,l\neq j}\varphi(q_l^{b_l})}{2^{k-3}}}\equiv  \pm 1 \mod q_j^{b_j},$$
\noindent because $\frac{\varphi(p_1^{a_1})}{4}\frac{\prod^{k-1}_{l=2,l\neq j}\varphi(q_l^{b_l})}{2^{k-3}}$ is odd. On the other hand,
$$2^{\frac{\varphi(n)}{2^k}}=(2^{\frac{\varphi(p_1^{a_1})}{4}})^
{\frac{\prod^{k-1}_{j=2}\varphi(q_j^{b_j})}{2^{k-2}}}\not \equiv +1 \textrm{ or } -1 \mod p_1^{a_1}.$$
\noindent Indeed, $p_1^{a_1} \notin (4+)\cup(4-)\cup(2+)$, because $p_1\equiv -3 \mod 8$; therefore,
$2^{\frac{\varphi(p_1^{a_1})}{4}}\not \equiv \pm 1 \mod p_1^{a_1}$,
and
$2^{\frac{\varphi(p_1^{a_1})}{2}}\not \equiv +1 \mod p_1^{a_1}$. Hence,

\noindent $2^{\frac{\varphi(n)}{2^k}}=(2^{\frac{\varphi(p_1^{a_1})}{4}})^
{\frac{\prod^{k-1}_{j=2}\varphi(q_j^{b_j})}{2^{k-2}}} \equiv +1\mod p_1^{a_1}$ only if $\frac{\prod^{k-1}_{j=2}\varphi(q_j^{b_j})}{2^{k-2}}$
is divisible by 4. But $\frac{\prod^{k-1}_{j=2}\varphi(q_j^{b_j})}{2^{k-2}}$ is odd, so $2^{\frac{\varphi(n)}{2^k}}\not \equiv +1 \mod 
p_1^{a_1}$.

Since $2^{\frac{\varphi(p_1^{a_1})}{2}} \equiv -1 \mod p_1^{a_1}$, we have $2^{\frac{\varphi(n)}{2^k}} \equiv -1$ only if 
$\frac{\prod^{k-1}_{j=2}\varphi(q_j^{b_j})}{2^{k-1}}$
is an odd integer. But $\frac{\prod^{k-1}_{j=2}\varphi(q_j^{b_j})}{2^{k-2}}$ is odd, which implies that 
$2^{k-1} \nmid \prod^{k-1}_{j=2}\varphi(q_j^{b_j})$; therefore, 
$2^{\frac{\varphi(n)}{2^k}}\not \equiv -1 \mod p_1^{a_1}$. 

So we have $2^{\frac{\varphi(n)}{2^k}}\equiv \pm 1 \mod q_j^{b_j}$ and $2^{\frac{\varphi(n)}{2^k}}\not \equiv \pm 1 \mod p_1^{a_1}$,
which implies that
$2^{\frac{\varphi(n)}{2^k}}\not \equiv \pm 1 \mod n$, i.e. $n \notin (2^k+)\cup(2^k-)$.

If $r=2$, then $n=p_1^{a_1}p_2^{a_2}\prod^{k-1}_{j=3}q_j^{b_j}$ and $\frac{\varphi(n)}{2^k}=
\frac{\varphi(p_1^{a_1})\varphi(p_2^{a_2})}{8}(\frac{\prod^{k-1}_{j=3}\varphi(q_j^{b_j})}{2^{k-3}})$. So
$$2^{\frac{\varphi(n)}{2^k}}=(2^{\frac{\varphi(p_1^{a_1})}{2}})^{\frac{\varphi(p_2^{a_2})}{4}
\frac{\prod^{k-1}_{j=3}\varphi(q_j^{b_j})}{2^{k-3}}}\equiv (-1)^{\frac{\varphi(p_2^{a_2})}{4}
\frac{\prod^{k-1}_{j=3}\varphi(q_j^{b_j})}{2^{k-3}}}\equiv -1 \mod p_1^{a_1},$$
\noindent because $\frac{\varphi(p_2^{a_2})}{4} \frac{\prod^{k-1}_{j=3}\varphi(q_j^{b_j})}{2^{k-3}}$ is odd. Also,
$$2^{\frac{\varphi(n)}{2^k}}=(2^{\frac{\varphi(p_2^{a_2})}{2}})^{\frac{\varphi(p_1^{a_1})}{4}
\frac{\prod^{k-1}_{j=3}\varphi(q_j^{b_j})}{2^{k-3}}}\equiv (-1)^{\frac{\varphi(p_1^{a_1})}{4}
\frac{\prod^{k-1}_{j=3}\varphi(q_j^{b_j})}{2^{k-3}}}\equiv -1 \mod p_2^{a_2}.$$
\noindent On the other hand,$$2^{\frac{\varphi(n)}{2^k}}=(2^{\frac{\varphi(q_j^{b_j})}{2}})^{\frac{\varphi(p_1^{a_1})\varphi(p_2^{a_2})}{8}
\frac{\prod^{k-1}_{l=3,l\neq j}\varphi(q_l^{b_l})}{2^{k-4}}}\equiv (\pm 1)^{\frac{\varphi(p_1^{a_1})\varphi(p_2^{a_2})}{8}
\frac{\prod^{k-1}_{l=3,l\neq j}\varphi(q_l^{b_l})}{2^{k-4}}}$$ 
$$\equiv +1 \mod q_j^{b_j},\textrm{ because } \frac{\varphi(p_1^{a_1})\varphi(p_2^{a_2})}{8}\textrm{ is even.}$$
\noindent Therefore, $2^{\frac{\varphi(n)}{2^k}}\equiv +1 \mod q_j^{b_j}$
and $2^{\frac{\varphi(n)}{2^k}}\equiv -1 \mod p_i^{a_i}$ which contradicts both $2^{\frac{\varphi(n)}{2^k}}\equiv +1 
\textrm{ and  } -1 \mod n$. 

Thus, $n \notin (2^k+)\cup(2^k-)$.
\\[10pt]
If $2 < r < k-1$, then 
$$2^{\frac{\varphi(n)}{2^k}}=(2^{\frac{\varphi(p_i^{a_i})}{2}})^
{\frac{\prod^{r}_{l=1,l\neq i}\varphi(p_l^{a_l})}{2^{r}}\frac{\prod^{k-1}_{j=r+1}\varphi(q_j^{b_j})}{2^{k-r-1}}}\equiv
(-1)^{\frac{\prod^{r}_{l=1,l\neq i}\varphi(p_l^{a_l})}{2^{r}}\frac{\prod^{k-1}_{j=r+1}\varphi(q_j^{b_j})}{2^{k-r-1}}}$$
$$\equiv +1 \mod p_i^{a_i}, \textrm{ because } \frac{\prod^{r}_{l=1,l\neq i}\varphi(p_l^{a_l})}{2^{r}}\textrm{ is even. }$$
$$\textrm{ Also, }2^{\frac{\varphi(n)}{2^k}}=(2^{\frac{\varphi(q_j^{b_j})}{2}})^
{\frac{\prod^{r}_{i=1}\varphi(p_i^{a_i})}{2^{r+1}}\frac{\prod^{k-1}_{l=r+1,l\neq j}\varphi(q_l^{b_l})}{2^{k-r-2}}}
\equiv (\pm 1)^{\frac{\prod^{r}_{i=1}\varphi(p_i^{a_i})}{2^{r+1}}\frac{\prod^{k-1}_{l=r+1,l\neq j}\varphi(q_l^{b_l})}{2^{k-r-2}}}$$
$$\equiv +1 \mod q_j^{b_j},\textrm{ because }\frac{\prod^{r}_{i=1}\varphi(p_i^{a_i})}{2^{r+1}}\textrm{ is even.}$$
\noindent So 
$2^{\frac{\varphi(n)}{2^k}}\equiv +1 \mod p_i^{a_i}$ and $2^{\frac{\varphi(n)}{2^k}}\equiv +1 \mod q_j^{b_j}$ for all $i$ and $j$.
Thus, $2^{\frac{\varphi(n)}{2^k}}\equiv +1 \mod n$, or $n\in (2^k+)$.
\hspace{\stretch{1}} $\Box$
\\[10pt]
We can easily see that Theorems \ref{th:7}, \ref{th:8}, \ref{th:9}, and \ref{th:10} 
cover all possible cases of $n$ with $k-1\geq 2$ different prime divisors. Therefore, the proof of Theorem 2 in the introduction
is complete.
\\[10pt]
{\bf Remark: }Using our theorems, it is easy to check whether $n\in (2^k+),\ n\in (2^k-)$, or $n\notin (2^k+)\cup(2^k-)$, in case when  $n$ has $k-1$ or more
distinct prime divisors. If $n$ has less than $k-1$ distinct prime divisors, then we can easily see that $n\notin (2^k+)\cup(2^k-)$,
if $n=\prod^r_{i=1}p_i^{a_i}$, with $0 < r < k-1$ and $p_i \equiv -1$ or 3 mod 8 for all $i$, because in this case $2^k \not | \varphi(n)$.
\section{The classes $(2\pm),(4\pm)$ and $(8\pm)$}
The classes $(2\pm)$ or the case $k=1$ is completely described in section \ref{sec:2}. 
Now let us study the case $k=2$ or the classes $(4\pm)$.

The following properties are proved in \cite{Ferm-Eul}, and we will see that, they can be obtained as special cases 
of the theorems in section \ref{sec:2}.
\begin{enumerate}
\item[P1.] Every odd number, divisible by more than two primes, belongs to the class $(4+)$.
\\[10pt]
The following 5 properties describe the case in which $n$ is divisible by two different primes $p$, $q$.
\item[P2.] If $p \equiv +1\mod 8,\textrm{ then } n \in (4+)$.
\item[P3.] If $p \equiv q \equiv -1\mod 8 \textrm{ or } p\equiv q\equiv -3\mod 8,\textrm{ then } n \in (4+)$. 
\item[P4.] If $p\equiv q\equiv +3\mod 8,\textrm{ then } n \in (4-)$.
\item[P5.] If $p \equiv +3\mod 8, q \equiv -1 \textrm{ or } -3\mod 8$, then $n \notin (4+)\cup(4-)$.
\item[P6.] If $p \equiv -3\mod 8, q \equiv -1 \mod 8$, then $n \notin (4+)\cup(4-)$.
\\[10pt]
The following property considers the case in which $n$ has only one prime divisor.
\item[P7.] If $n = p^a$ belongs to the class $(4+)\cup(4-)$, then $p \equiv 1\mod 8$.
\\[10pt]
The following two properties consider every odd number $n$.
\item[P8.] If $n \equiv 5 \mod 8$, then $n$ does not belong to the class $(4-)$.
\item[P9.] If $n \equiv 7 \mod 8$, then $n$ does not belong to the class $(4-)$.
\end{enumerate}
\begin{itemize}
\item It is clear that if $k=2$, then P1, P2, P3, and P4 are special cases from theorems \ref{th:1}, \ref{th:2}, \ref{th:3}, and
\ref{th:4} respectively.

\item In P5 we have two cases:

Case 1: If $p \equiv +3\mod 8, q \equiv -1\mod 8$, which is exactly theorem \ref{th:6} with $k=2$.

Case 2: If $p \equiv +3\mod 8, q \equiv -3\mod 8$, which is exactly the case in theorem \ref{th:5} with $k=2$ and $r=1$.

\item The case in P6 is the case in theorem \ref{th:5} with $k=2$ and $r=1$.

\item In P7, we have just one direction and Lemma \ref{th:3.1} complete the second direction.
\end{itemize}
\begin{itemize}

\item In P8 and P9, $n$ has three cases. Case 1: If $n$ has more than two distinct prime divisors, then by theorem \ref{th:2}, 
$n$ belongs to $(4+)$ not to $(4-)$.

Case 2: If $n$ has two prime divisors, then by theorem \ref{th:4}, $n$ belongs to $(4-)$ only if all the prime divisors of $n$
are congruent to 3 mod 8. Since $3^2 \equiv 1 \mod 8,\textrm{ we have } n \in (4-) \textrm{ only if } n \equiv 3 \textrm{ or } 1 \mod 8$.
So $n\notin (4-)$ if $n\equiv 5$ or 7 mod 8.

Case 3: If $n$ has only one prime divisor, then by Lemma \ref{th:3.1} $n \notin (4+)\cup(4-).$
\end{itemize}

\noindent {\bf The classes $(8\pm)$:}
\\[10pt]
If $k=3$, theorem \ref{th:1} state the following:

\noindent Every odd number, divisible by more than three primes, belongs to the class $(8+)$.

So now consider the case $n=p^aq^bs^c, a>0, b>0, c>0$. In \cite{Ergodic-arithmatic}, Arnold denoted by I, II, III, IV the sets
of the primes, congruent to 1, 3, 5 and to 7 mod 8 respectively, and proved the following properties:

\begin{enumerate}
\item[P1.] The triple product $n$ belongs to the class $(8-)$, whenever the primes $(p, q, s)$ 
do all belong to the kind II.

\item[P2.] The triple product $n$ belongs to the class $(8+)$, provided that the
three primes $p, q$ and $s$ do belong to one of the following 7 triples of kinds:

(I, I, X), (I, III, X), (I, II, II), (I, II, IV), (I, IV, IV), (IV, IV, IV), 

(III, III, X).

\item[P3.] No triple product $n$ belongs to the class $(8+)$ or $(8-)$, provided
that the primes $(p, q, s)$ kinds triple is one of the 5 triples:

(III, II, II), (III, II, IV), (III, IV, IV), (II, IV, II), (II, IV, IV) .
 
\end{enumerate}
Note that the triples in P1, P2, and P3 cover (taking the permutations into account) all the 
64 possible ordered lists of the 3 kinds of $p$, $q$, and $s$. To verify it, see \cite{Ergodic-arithmatic}.
\\[10pt]
Now consider the case if $n$ has exactly $k-1=2$ distinct prime divisors. Let $n=p^aq^b, a>0,b>0$. In \cite{Ergodic-arithmatic}
Arnold proved the following properties:
\begin{enumerate}

\item[P4.] If $p \in$ III and $q \in$ III, then $n\in (8-)$.

\item[P5.] If $p \in$ I and $q \in$ I $\cup$ III, then $n \in (8+)$.

\item[P6.] If $p \in$ (I $\cap (4+))$ and $q \in$ (II $\cup$ IV), then $n \in (8+)$.

\item[P7.] If $p \in$ (I $\cap (4-))$ and $q \in$ (II $\cup$ IV), then $n \notin (8+)\cup(8-)$.

\item[P8.] If $p\in$ (II $\cup$ IV) and $q \in$ (II $\cup$ IV), then $n \notin (8+)\cup(8-)$.
 
\item[P9.] If $p\in (4+) \cup (4-)$ and $q\in (4+) \cup (4-)$, then $pq \in (8+)$.
 
\end{enumerate}

We will see now that these properties are special cases of the theorems in section \ref{sec:2}.
\begin{itemize}
\item It is clear that P1 is exactly theorem \ref{th:4} with $k=3$.

\item In P2, the first 5 triples (I, I, X), (I, III, X), (I, II, II), (I, II, IV), (I, IV, IV) have at least one prime divisor
of kind I, that is they satisfy the condition of theorem \ref{th:2}. So they are special cases of it.

The triple (IV, IV, IV) is exactly theorem \ref{th:3} with $k=3$.

The triple (III, III, X) has 3 cases. Case 1: the triple (III, III, I) clearly belongs to theorem \ref{th:2}.

Case 2: the triples (III, III, II) and (III, III, IV)  have two prime divisors of kind 
III and only one divisor of kind II or IV, so they belong to theorem \ref{th:5} with $k=3, r= 1<k-1$. 

Case 3: the triple (III, III, III) is exactly theorem \ref{th:3} with $k=3$.

\item In P3, the first 3 triples (III, II, II), (III, II, IV), (III, IV, IV) have one prime divisor from kind III and the others
from kind II or IV, which means that they are special cases from theorem \ref{th:5} with $k=3, r=k-1=2$.

The triples (II, IV, II), (II, IV, IV) have only prime divisors from kinds II or IV, so they belong to 
theorem \ref{th:6} with $k=3$.

\item It is clear that the properties P4-P9 can be obtained from Theorems \ref{th:7}-\ref{th:10} and Lemma 1. 
\end{itemize}



\begin{thebibliography}{10}

\bibitem{Ferm-Eul} \textsc{Arnold, V.~I.}, \emph{Fermat-Euler Dynamical Systems and the Statistics of Arithmetics of Geometric Progressions}, 
 Funktsional. Anal. i Prilozhen.  {\bf 37}  (2003),  no. 1, 1--18, 95;  translation in  Funct. Anal. Appl. {\bf 37}  (2003),  no. 1, 1--15.

\bibitem{Ergodic-arithmatic} \textsc{Arnold, V.~I.}, \emph{Ergodic and Arithmetical properties of geometrical progression's dynamics and its orbits},
 Mosc. Math. J.  {\bf 5}  (2005),  no. 1, 5--22.

\bibitem{Friendly} \textsc{Silverman, Joseph H.}, \emph{ A Friendly Introduction to Number Theory}, Prentice Hall, {\bf 1996}.
\end{thebibliography}
\end{document}